\documentclass[trsc,nonblindrev]{informs3-neutral}

\OneAndAHalfSpacedXI 

\usepackage{float}
\usepackage{booktabs}
\usepackage[flushleft, para]{threeparttable}
\usepackage{algorithm} 
\usepackage{mathtools}
\usepackage{algpseudocode}
\algnewcommand{\NoEndIf}[1]{\State\algorithmicif\ #1\ \algorithmicthen\ }
\usepackage{enumitem}
\usepackage{comment}
\usepackage{commands}
\usepackage[nolist]{acronym}
\usepackage{placeins}
\usepackage{hyperref}

\newtheorem{prop}{Proposition}

\newtheorem{result}{Result}

\makeatletter
\def\tcb@cnt@exampleautorefname{Result}
\makeatother

\newcommand*{\QEDB}{\null\nobreak\hfill\ensuremath{\square}}%
\makeatletter
\renewcommand\paragraph{\@startsection{paragraph}{4}{\z@}%
                                    {3.25ex \@plus1ex \@minus.2ex}%
                                    {-1em}%
                                    {\normalfont\normalsize\bfseries}}
\makeatother

 \usepackage {tikz}

 \usepackage[colorinlistoftodos,textsize=scriptsize]{todonotes}

\usepackage{natbib}
 \bibpunct[, ]{(}{)}{,}{a}{}{,}%
 \def\bibfont{\small}%
\TheoremsNumberedThrough     

\EquationsNumberedThrough    

\begin{document}



\RUNTITLE{CO and ML for Dynamic Inventory Routing}

\TITLE{Combinatorial Optimization and Machine Learning for Dynamic Inventory Routing}

\ARTICLEAUTHORS{%
\AUTHOR{Toni Greif}
\AFF{Julius-Maximilians-Universität Würzburg, Germany
\EMAIL{toni.greif@uni-wuerzburg.de}}
\AUTHOR{Louis Bouvier}
\AFF{CERMICS, \'Ecole des Ponts, Marne La Vallée, France, \EMAIL{louis.bouvier@enpc.fr}}
\AUTHOR{Christoph M. Flath}
\AFF{Julius-Maximilians-Universität Würzburg, Germany
\EMAIL{christoph.flath@uni-wuerzburg.de}}
\AUTHOR{Axel Parmentier}
\AFF{CERMICS, \'Ecole des Ponts, Marne La Vallée, France, \EMAIL{axel.parmentier@enpc.fr}}
\AUTHOR{Sonja U. K. Rohmer}
\AFF{Department of Logistics and Operations Management, HEC Montréal, Canada, \EMAIL{sonja.rohmer@hec.ca}}
\AUTHOR{Thibaut Vidal}
\AFF{CIRRELT \& SCALE-AI Chair in Data-Driven Supply Chains, MAGI, Polytechnique Montréal, \EMAIL{thibaut.vidal@polymtl.ca}} 

} 

\ABSTRACT{%
We introduce a combinatorial optimization-enriched machine learning pipeline and a novel learning paradigm to solve inventory routing problems with stochastic demand and dynamic inventory updates. After each inventory update, our approach reduces replenishment and routing decisions to an optimal solution of a capacitated prize-collecting traveling salesman problem for which well-established algorithms exist. Discovering good prize parametrizations is non-trivial; therefore, we have developed a machine learning approach. We evaluate the performance of our pipeline in settings with steady-state and more complex demand patterns. Compared to previous works, the policy generated by our algorithm leads to significant cost savings, achieves lower inference time, and can even leverage contextual information.
}


\maketitle

%


\begin{acronym}
\acro{IRP}{Inventory Routing Problem}
\acro{SIRP}{Stochastic Inventory Routing Problem}
\acro{DSIRP}{Dynamic and Stochastic Inventory Routing Problem}
\acro{SAA}{Sample Average Approximation}
\acro{CO}{Combinatorial Optimization}
\acro{ML}{Machine Learning}
\acro{MILP}{Mixed Integer Linear Program}
\acro{RL}{Reinforcement Learning}
\acro{MDP}{Markov Decision Processes}
\acro{TSP}{Traveling Salesman Problem}
\acro{CPCTSP}{Capacitated Prize Collecting TSP}
\acro{DVRP}{Dynamic Vehicle Routing Problem}
\acro{PINN}{Physics-informed Neural Network}
\end{acronym}

\section{Introduction}\label{dsirp:sec:intro}
The \ac{IRP} is a challenging \ac{CO} problem which simultaneously considers routing, inventory holding and stock-out costs. This problem is of high relevance for various real-world applications in transportation, logistics, and supply chain management.
There has been extensive research on the deterministic \ac{IRP}, resulting in various solution methods and heuristics. However, inventory management applications are typically characterized by different uncertainty sources which in turn affect different components of the \ac{IRP}. 
In this study, we focus on uncertainty within the customer demands, assuming that the demands of the current period only become known once the replenishment decisions have been made. As such, we consider both the stochasticity of demand as well as the dynamic nature of the problem requiring continuous inventory updates at the customer locations. Consequently, the problem state, which includes the inventory levels, historical demand observations, and possibly predictive contextual information for each customer, is consistently updated after each period. This combination of stochastic customer demands and inventory updates leads to the \ac{DSIRP}, which may lead to stock-outs even under optimal decision-making \citep{coelho2014heuristics}.

The goal of this paper is to introduce a learning-based policy for the \ac{DSIRP}, and benchmark it against classic rolling-horizon policies based on single or multiple demand scenarios.
Single-scenario policies are well established and easy to apply, but typically lack performance when demand patterns get more complex.
Multi-scenario policies overcome these performance issues.
However, their computational complexity increases rapidly with the number of scenarios and customers considered, which makes them impractical even for small instances. Our learning-based approach aims to combine the strengths of \ac{CO} and \ac{ML} in a hybrid pipeline to achieve optimal decision-making under uncertainty.

As a result, our research represents a pioneering effort in addressing the \ac{DSIRP} using advanced algorithms. Vis-a-vis traditional solution approaches, our learning-based pipeline offers strong performance across a large number of problem configurations while retaining short inference times. This combination enables robust real-time decision-making in complex \ac{DSIRP} settings. The strength of our pipeline lies in its resilience and adaptability. This is because it operates without being affected by the underlying demand distribution but instead relies solely on historical samples. Moreover, our approach has the advantage of effectively leveraging contextual information, which sets it apart from other algorithms that either only have limited capabilities to incorporate such information or do not scale well in terms of computational complexity. We provide the instances used in our research to support further studies on contextual demand in the \ac{DSIRP} literature. In addition, while the potential of integrating hybrid \ac{ML} and \ac{CO} pipelines has been highlighted in simplified settings, our work is among the first attempts, together with a few notable exceptions, to address a complex \ac{IRP} using such methods. By thoroughly analyzing the strengths and limitations of our approach, we offer guidance to both practitioners and researchers operating in this field.

This paper provides insights into the performance, robustness, and real-time decision-making capabilities of the presented methods, shedding light on the inherent trade-offs of selecting an appropriate approach to effectively solve the \ac{DSIRP} for real-world scenarios. Focusing on a simplified setting, our numerical experiments allow us furthermore to evaluate the effectiveness of the most complex approaches in finding optimal solutions. Moreover, by isolating the heuristic solutions from other confounding factors, we are able to analyze and compare the relative strengths of these approaches.

We refer to our git repository \url{https://github.com/tonigreif/InferOpt_DSIRP.git} for all the material necessary to reproduce the results outlined in this paper.

\section{Related Work}\label{sec:related_work}
Contributing to two research directions, we first review the \ac{IRP} and  \ac{DVRP} literature before investigating the potential of \ac{ML}-based approaches to address stochastic and dynamic problems. In the latter, we focus on the concept of \textit{\ac{CO}-enriched \ac{ML}}, which presents a hybrid pipeline where \ac{ML} models learn from many similar \ac{CO} problems.

\subsection{Inventory Routing}\label{sec:SIRP}
The \ac{IRP} is a well-known problem for which many extensions and modifications have been studied with the aim of representing a large variety of different real-world situations \citep{coelho2014thirty}. This has resulted in an abundance of solution approaches for the deterministic \ac{IRP}, including exact methods such as branch-and-cut \citep{adulyasak2014formulations, archetti2007branch,avella2018single, bertazzi2019matheuristic, guimaraes2019two, coelho2013branch, coelho2013exact, manousakis2021improved} and branch-and-price-and-cut \citep{desaulniers2016branch}, as well as heuristics and metaheuristics \citep{adulyasak2014formulations, archetti2012hybrid, archetti2017matheuristic,bertazzi2019matheuristic, chitsaz2019unified,guimaraes2019two,bouvier2022solving}.

Given the importance of uncertainty in real-world settings, there is an increasing need to reconsider the deterministic nature of \ac{IRP}s and adequately account for stochasticity when solving such problems. The review of \citet{soeffker2022stochastic} identifies in this context three key problem dimensions that may be affected by uncertainty: i) demand (requests, quantities, service time, etc.), ii) resources (vehicle availability, driver availability, range, etc.), and iii) the environment (travel times, fees, road closures, etc.). Arising in almost all practical settings, demand uncertainty is one of the most common sources of uncertainty considered in the scientific literature.
Accordingly, this research will also focus on incorporating demand uncertainty and the most relevant studies in this area. 
A common approach to deal with demand uncertainty without information on the underlying probability distribution is to estimate customer demands based on historical samples.
A straightforward strategy for addressing these \ac{CO} problems is selecting multiple (demand) scenarios from these historical samples and formulating the \ac{MILP}. This approach involves introducing what we refer to as first-period linkage constraints. These constraints play a crucial role during the problem-solving phase, ensuring the uniqueness of the first-period routing decision across all scenarios. To handle the considerable increase in size and complexity, consensus optimization frameworks avoid using first-period linkage constraints. Instead, they decompose the \textit{global consensus problem} into independent subproblems, iterating to achieve consensus on the solution \citep{nedic201010,boyd2011distributed, bertsekas2015parallel}. Under this paradigm, multiple \ac{IRP} scenarios are solved independently, and consensus strategies like progressive hedging are employed to merge divergent routing solutions into a single decision \citep{hvattum2009using}.

The least complex option involves consolidating multiple scenarios into a single scenario before solving the problem, as demonstrated in the case of the \ac{DSIRP} by \citet{coelho2014heuristics}. This option is often employed when the single-scenario problem is intricate in itself. As this approach may lack effectiveness, the study by \citet{brinkmann2019dynamic} proposes simulation techniques that involve multiple scenarios. In addition, different heuristics were proposed to solve this type of problem; examples of some notable methods, in this context, include ant colony optimization algorithms \citep{huang2010modified}, variable neighborhood search metaheuristic hybridized with simulation \citep{gruler2018combining}, a hybrid rollout algorithm \citep{bertazzi2013stochastic}, and benders decomposition \citep{li2022two}. 

In addition to stochastic considerations, uncertainty may add dynamic considerations to a problem as the initial stochastic information is updated or revealed over time. 
Integrating dynamic aspects in \ac{IRP}s is relatively new and has not been widely explored within the literature. For this reason we present a more general overview of dynamic vehicle routing problems. Examples of these problems in different routing contexts are provided by \citet{ritzinger2016survey, gendreau201650th, psaraftis2016dynamic, oyola2018stochastic}. 
A common strategy for solving these problems is the rolling-horizon approach, where a static solution is continually reoptimized as new information becomes available \citep{berbeglia2010dynamic, wong2006optimal, bouvier_dynamic_2023}.

\subsection{Decision-Making with Machine Learning}
In recent years, the use of \ac{ML} techniques for decision-making has received growing interest in the scientific literature. Several studies highlight the potential of such techniques in addressing \ac{DVRP}s \citep{joe2020deep, james2019online}.
\citet{HILDEBRANDT2023106071} emphasize, in particular, the potential of \ac{RL} for two compelling reasons: First, the ability to model dynamic \ac{CO} problems as \ac{MDP} \citep{powell2019unified}. Second, offline training of learning architectures to devise complex decision policies with short inference times.
However, \citet{HILDEBRANDT2023106071} also point out that \ac{ML} techniques face challenges when dealing with complex, combinatorial action spaces.
In such scenarios, the search for the action space is better performed by available \ac{CO} solvers. 
Hence, there is untapped potential in combining \ac{ML} with \ac{CO} within a hybrid pipeline. Such an approach can leverage the strengths of both paradigms and holds promise for developing combinatorial solutions while accounting for uncertainties.

These pipelines, known as \textit{\ac{CO}-enriched \ac{ML}} pipelines, have demonstrated remarkable success in solving problems that were challenging for classic \ac{CO} approaches. More specifically, they have shown excellent performance in addressing problems like single machine scheduling with release dates \citep{parmentier2021learning} and various multi-stage optimization problems. Examples of successful multi-stage optimization problems include the parameterization of the two-stage stochastic minimum spanning tree problem \citep{dalle2022learning}, the stochastic vehicle scheduling problem \citep{parmentier2021learningA, parmentier2022learning}, and even the dispatching policy in a dynamic autonomous mobility-on-demand system \citep{jungel2023learningbased}. Additionally, this method has proven effective in solving high-dimensional multi-stage stochastic optimization problems like the \ac{DVRP} with time windows \citep{baty2023combinatorial}.

\section{Problem Description}\label{dsirp:sec:problemDescription}
This paper considers the \ac{DSIRP}, a dynamic variant of the \ac{IRP} where information on customer demands and inventory updates are revealed progressively.
Customer demands are assumed to be random variables and may be interpreted as a retailer's total order for a given period.
\begin{figure}[htb]
     \centering
     \begin{tikzpicture}
       \tikzset{
         box/.style = {draw, align=center, rectangle, minimum width = 2.5cm, minimum height = 0.5cm},
         arrow/.style = {->, thick, >=stealth}
       }
      \draw (0,0) -- (1.5,0) node[right] {...};
      \draw (14,0) -- (15.5,0) node[left] at (14,0) {...};
    
      \draw (0,-0.2) -- (0,0.2) node[below=5mm ] {$t=0$};
      \foreach \timestep/\label in {1/1, 14.5/{T-1}, 15.5/T} {
        \draw (\timestep,-0.2) -- (\timestep,0.2) node[below=5mm] {$\label$};
      }
       \draw[arrow] (2.2, 0) -- node[above] {{Inventories}} (5, 0);
       \draw[arrow] (2.2, 0) -- node[below] {{$(I^t_i)_{i\in\mathcal{V}_c}$}} (5, 0);
       \node[box, anchor=west] at (5, 0) {Decisions \\ $(u^t, z^{t})$};
       \draw[arrow, anchor=west] (7.5, 0) -- node[above] {} (8, 0);
       \node[box, anchor=west] at (8, 0) {Demand\\ $d(\xi^t)$};
       \draw[arrow, anchor=west] (10.5, 0) -- node[above] {Inventories} (13.3, 0);
       \draw[arrow, anchor=west] (10.5, 0) -- node[below] {$(I^{t+1}_i)_{i\in\mathcal{V}_c}$} (13.3, 0);
     \end{tikzpicture}
     \caption{Process of making decisions, revealing demand, and updating inventories over time}
     \label{fig:decision_process}
 \end{figure}
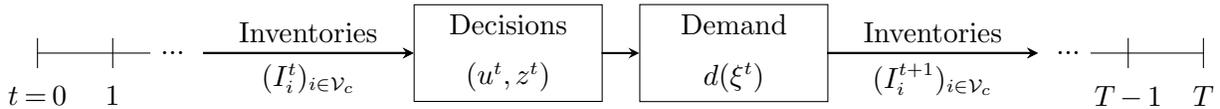

The problem is defined over a given \textit{planning horizon}, at the beginning of which each customer has an initial inventory. Inventory positions incur a customer-specific holding cost per period and unit. In the same fashion, inventory shortfalls result in customer-specific stock-out cost per unit of unsatisfied demand. We do not consider backlogging of unfulfilled demand. To determine the replenishment decisions for each customer we use an order-up-to policy which is widely used in \ac{IRP} settings. Under this policy delivered quantities completely replenish the inventory capacity at the customer locations.
Production capacity is unconstrained we do not account for inventory holding costs at the supplier. As common in the \ac{IRP} literature, the distribution part of our problem consists of a single vehicle with a given capacity, which can perform one route per period. Each route starts at the supplier and covers a subset of customer locations, while traveling between locations incurs a travel cost. The problem's overall objective is to minimize the total routing cost and the costs of inventory holding and stock-outs at the customer locations for the total \textit{planning horizon}.

The problem can be defined on an undirected graph $\mathcal{G} = (\mathcal{V}, \mathcal{E})$, where $\mathcal{V}=\{0,\dots, v\}$ presents the set of vertices and $\mathcal{E}=\{(i,j):i,j \in \mathcal{V}, i < j\}$ the set of feasible edges. Vertex $0$ represents, in this context, the depot at which the supplier is located, and the vertices $\mathcal{V}_c=\mathcal{V}$\textbackslash$\{0\}$ correspond to the customer locations. 

\subsection{Markov Decision Process}
\label{subsec:MDP}

To provide the necessary theoretical foundation for our learning-based pipeline, we formulate the \ac{DSIRP} as a \ac{MDP}. For this purpose, we define the vector $I^{t} \in \mathbb{R}_{\geq 0}^{|\mathcal{V}_c|}$, which represents the inventory levels of the customers at the beginning of period $t$, not taking into account the forthcoming realized demand or the replenishment activities. The inventory level of a specific customer $i$ is then denoted by the element $I^{t}_i$.
Furthermore, we denote the state before making a decision at period $t$ as $x^t$. In the absence of additional contextual information, $x^t$ corresponds to the customers' inventory levels at the beginning of the current period, along with the set of the most recent historical observations $\hat{\mathcal{D}^t}=\{d^{t-1}, \dots ,d^{t-50}\}$ for each customer $i\in\mathcal{V}_c$. However, when contextual information $\varPhi^t$ regarding the environment for period $t$ is available, the state is augmented to $x^t=(I^t, \hat{\mathcal{D}^t}, \varPhi^t, \varPhi^{t+1}, \dots)$ to encompass this supplementary information for the subsequent periods, with $\varPhi_i^{t}$ representing the context for customer $i$.
We proceed from the premise that the context $\varPhi^t$ and the exogenous noise $\xi^t$ capturing the uncertainty in demands $d^t = d(\xi^t)$ are not only independent of each other, but also of all decisions made. This basic assumption ensures that we are in a Markovian setting.

We then denote by $\mathcal{T}$ the set of tours (without subtours) that start and end in vertex $0$. A tour can be encoded by a vector ($u^t_{ij})_{(i,j) \in \mathcal{E}}$, where $u^t_{ij} = 1$ if edge $(i,j)$ belongs to the tour and $0$ otherwise. 
To encode if customer $i$ belongs to the tour or not, we define
\begin{equation}\label{eq:func_g}
z^t = g(u^{t}) = \frac{1}{2}\left(\sum_{j\in \mathcal{V}, i>j}u_{ji}^t  +  \sum_{j\in \mathcal{V}, i<j}u_{ij}^t \right)_{i\in\mathcal{V}_c}.
\end{equation}
When the process occupies state $x^t$, the set of feasible decisions is then given by
\begin{equation}
\mathcal{U}(x^t) = \left\{
  u^t \in \mathcal{T} \;\middle|\;
  \sum_{i\in \mathcal{V}_c} \left(C_i - I^{t}_i\right) {z_i^{t}} \leq B
\right\},
\label{eq:decisionSet}
\end{equation}
indicating that we can replenish all subsets of customers whose combined replenishment quantities, determined by the difference between inventory capacity $C_i$ and current inventory $I_i^{t}$, are smaller than the vehicle capacity $B$.
Furthermore, the calculation of routing costs is defined by
\begin{equation}\label{eq:func_h}
h(u^{t})=\sum_{(i,j) \in \mathcal{E}} \gamma_{ij} u_{ij}^t,
\end{equation}
where $\gamma_{ij}$ are the usage costs of edge $(i,j) \in \mathcal{E}$.
Given state $x^t$ with inventory levels $I^t$, each decision is associated with a (negative) reward equal to its respective holding, stock-out, and routing costs:
\begin{equation}
    \label{eq:reward}
        \tilde r(x^{t}, u^t, \xi^t) = \sum_{i\in\mathcal{V}_c} \kappa_i \left\{I_i^{t} (1-z_i^{t}) + C_i z_i^{t} - d(\xi^t)_i\right\}^+ - \rho \kappa_i \left\{I_i^{t} (1-z_i^{t}) + C_i z_i^{t} - d(\xi^t)_i\right\}^- + h(u^{t}),
\end{equation}
where $\kappa_i$ are the holding costs per unit and the shortage penalty $\rho$ indicates the multiplier by which the holding costs influence the stock-out costs.
Due to the absence of recourse actions in case a customer runs out of stock, the deterministic transition to the next state $x^{t+1}$ corresponds to updating the customer inventory levels
\begin{equation}\label{eq:inventory_transition}
    I_i^{t+1} = \left\{I_i^{t} (1-{z_i^{t}}) + C_i {z_i^{t}} - d(\xi^t)_i\right\}^+,
\end{equation}
as well as the historical demand observations and, if available, also the contextual information, resulting in 
\begin{equation}\label{eq:state_transition_contextual}
x^{t+1} = (\underbrace{\vphantom{\varPhi^{t+1}, \varPhi^{t+2}, \dots} I^{t+1}}_{\text{inventories}},
  \underbrace{\vphantom{\varPhi^{t+1}, \varPhi^{t+2}, \dots} \hat{\mathcal{D}}^{t+1}}_{\text{historical demand}},
  \underbrace{\varPhi^{t+1}, \varPhi^{t+2}, \dots}_{\text{contextual information}}
).
\end{equation}

A deterministic policy $\delta$ then maps a state $x^t$ to a feasible decision $u^t \in \mathcal{U}(x^t)$. We seek a policy in the set of all Markovian deterministic policies $\Delta$ that minimizes the expected total (negative) reward conditional on the initial state $x^0$

\begin{equation}\label{eq:MDP_policy}
    \mathop{\arg \min}_{\delta\in\Delta}{\mathbb{E} \left[\sum_{t=0}^{T-1} \tilde r\big(x^t, \delta\left(x^t\right), \xi^t\big) \middle| x^0 \right]}.
\end{equation}

\section{Policy Encoded as Machine Learning Pipeline}
\label{sec:pipeline_policy}

The challenge of the \ac{DSIRP} arises due to the combinatorially vast nature of \textit{both} the state space $\calX$ and the decision space $\calU$. 
While classic stochastic optimization policies can handle large $\calU$ when $\calX$ is small, and reinforcement learning can handle large $\calX$ in the context of small $\calU$, settings with both large $\calX$ and large $\calU$ are still poorly addressed.
As a result, solution methods for \ac{DVRP}s often rely on rolling-horizon policies as discussed in \autoref{sec:related_work}.
Although these policies are practical and can account for uncertainty by considering multiple historical samples through voting or consensus strategies, computational limitations commonly restrict the number of scenarios that can be considered, reducing their effectiveness. 

\begin{figure}[htb]
    \centering
    \begin{tikzpicture}
      \tikzset{
        box/.style = {draw, align=center, rectangle, minimum width = 4cm, minimum height = 0.5cm},
        arrow/.style = {->, thick, >=stealth}
      }
      \draw[arrow] (0, 0) -- node[above] {State} (2, 0);
      \draw[arrow] (0, 0) -- node[below] {$x^t$} (2, 0);
      \node[box, anchor=west] at (2, 0) {Statistical model\\ $\varphi_w$};
      \node[above] at (4, 0.7) {\textbf{ML layer}};
      \draw[arrow, anchor=west] (6, 0) -- node[above] {Prizes} (9, 0);
      \draw[arrow, anchor=west] (6, 0) -- node[below] {$\theta = \varphi_w\left(x^t\right)$} (9, 0);
      \node[box, anchor=west] at (9, 0) {Oracle\\ $o$};
      \node[above] at (11, 0.7) {\textbf{CO layer}};
      \draw[arrow, anchor=west] (13, 0) -- node[above] {Decision} (16, 0);
      \draw[arrow, anchor=west] (13, 0) -- node[below] {$u^t=o(\theta)$} (16, 0);
    \end{tikzpicture}
    \caption{Our CO-enriched ML pipeline.}
    \label{fig:pipeline}
\end{figure}
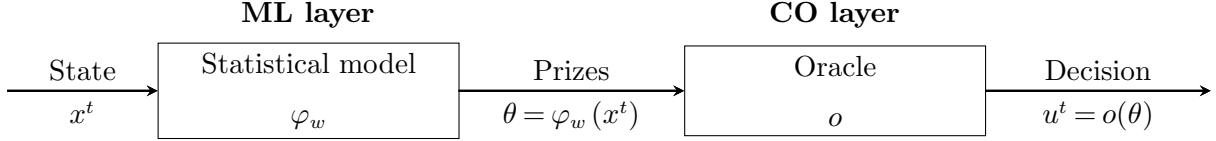

To overcome these limitations, we propose a hybrid pipeline for encoding policies. The pipeline chains a statistical model $\varphi_w$ with a standard \ac{CO} problem for which well-established algorithms exist. We refer to this solution algorithm as oracle $o$. The oracle must share the same set of feasible solutions $\mathcal{U}(x^t)$ and be parameterizable. With fixed parameters $\theta$ in each state $x^t$, the oracle's decision $u^t$ becomes a deterministic policy. This pipeline (\autoref{fig:pipeline}) is a \textit{\ac{CO}-enriched \ac{ML}} pipeline. It encodes a family of policies parametrized by $w$, which denotes the parameters of the statistical model $\varphi_w$ to predict $\theta$ given a state $x^t$.
An important task in this context is choosing the parameters $w$ such that the pipeline leads to a practically efficient policy. For this purpose, we introduce the learning algorithm, which is described in \autoref{sec:learning}.

In the following subsections, we will first introduce and explain our choice for the oracle $o$, formulating it as a \ac{CPCTSP}, and then our statistical model $\varphi_w$, a \ac{PINN}.

\subsection{Combinatorial Optimization Layer}
\label{subsec:CPCTSP}
The \ac{CPCTSP} is a deterministic vehicle routing problem variant in which customer deliveries are optional but rewarded by a prize. It aims to select a subset of the customers along with good delivery routes to maximize the sum of rewards minus transportation costs (distance). There exists close connections between this problem and the \ac{DSIRP} since the set of feasible decisions $\mathcal{U}(x^t)$, at each state $x^t$, consists of tours that visit a subset of customers whose combined replenishment quantities do not exceed the vehicle capacity, as stated in Equation~\eqref{eq:decisionSet}. This set aligns with the feasible solutions of a \ac{CPCTSP}, parameterized by prizes $\theta$ representing how desirable deliveries are to the different locations at the considered period. Moreover, solving a single-period \ac{CPCTSP} is less complex than a multi-period \ac{IRP} and benefits from well-established algorithms. As such, we define our problem as follows.
Let $\mathcal{G} = (\mathcal{V}, \mathcal{E})$ be the undirected graph defined in \autoref{dsirp:sec:problemDescription}. Each vertex $i \in \mathcal{V}_c$ has an associated prize $\theta_i$. For each edge $(i,j)$ in $\mathcal{E}$, we have a travel cost $\gamma_{ij} \geq 0$.
The task is then to determine the customers to be visited, by constructing a tour that begins at a depot. The goal is to maximize profit, calculated as the sum of the collected prizes minus the travel costs incurred on the traversed edges.
Using Equation~\eqref{eq:func_g}, the resulting \ac{CPCTSP} can be formulated as 
\begin{equation}
    \label{eq:CPCTSP}
    \max_{u \in \mathcal{U}} \underbrace{\sum_{i \in \mathcal{V}_c}\theta_i {z_i^t} }_{\text{prizes collected}} - \underbrace{\sum_{(i,j) \in \mathcal{E}} \gamma_{ij} u_{ij}^t}_{\text{routing costs}}.
\end{equation}
It is important to note that the routing costs are well-defined within our \ac{DSIRP} context, but the prizes are not. However, we demonstrate that there always exists a prize vector $\theta$ such that the optimal solution of the \ac{CPCTSP} aligns with the optimal decision for the \ac{DSIRP}.
\begin{prop}
For every state $x^t$ there exists a prize vector $\theta \in \mathbb{R}^{|\mathcal{V}_c|}$ such that any optimal solution of the \ac{CPCTSP} \eqref{eq:CPCTSP} corresponds to an optimal decision in terms of the expected total (negative) reward of the \ac{MDP} \eqref{eq:MDP_policy}.
\end{prop}
\begin{proof}{Proof.}
Since the planning horizon is finite and the set of feasible decisions at each period is also finite, there exists an optimal decision $u^*$ for $x^t$. Let $\bar{\mathcal{V}}_c$ be the subset of customers of $\mathcal{V}_c$ that are replenished in $u^*$. Then any solution $u$ which has lower or equally lower routing costs and covers $\bar{\mathcal{V}}_c$ exactly is also optimal. This follows from the Bellman equation since the routes have no impact on the evolution of the state. We can construct $u$ by solving a \ac{CPCTSP} on $\mathcal{V}_c$ with prizes denoted by
\begin{equation}
    \theta_i = \begin{cases}
+M &\text{if } i\in \bar{\mathcal{V}}_c, \\
-M & \text{otherwise,}
\end{cases}
\end{equation}
where $M=|\mathcal{V}_c| \cdot \max_{(i,j)\in \mathcal{E}} \gamma_{ij}$ is a large constant. The corresponding \ac{CPCTSP} solution $u$ clearly covers $\bar{\mathcal{V}}_c$ and has at most the routing costs of $u^*$.
\QEDB 
\end{proof}

\subsection{Machine Learning Layer}
\label{subsec:Stat-Model}

The goal of the \ac{ML} layer is to predict for any state $x^t$ a prize vector $\theta$ in such a way that optimal solutions of the resulting \ac{CPCTSP} align with good solutions of the \ac{DSIRP} for this state.

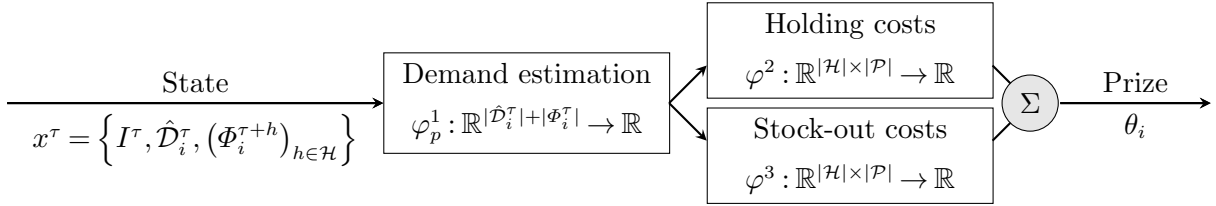
\begin{figure}[htb]
    \centering
    \begin{tikzpicture}[
        sum/.style = {draw, circle, minimum size=0.7cm, fill=black!10},
        box/.style = {draw, align=center, rectangle, minimum width=3.8cm, minimum height=0.5cm},
        arrow/.style = {->, thick, >=stealth}
    ]
      \draw[arrow] (0, 0) -- node[above] {State} (5, 0);
      \draw[arrow] (0, 0) -- node[below] {$x^\tau = \left\{I^\tau, \hat{\mathcal{D}}_i^\tau, \left(\varPhi_i^{\tau+h}\right)_{h\in\mathcal{H}}\right\}$} (5, 0);
      \node[box, anchor=west] at (5, 0) {Demand estimation\\ $\varphi^1_p: \mathbb{R}^{|\hat{\mathcal{D}}_i^\tau| + |\varPhi_i^{\tau}|} \rightarrow \mathbb{R}$};
      \draw[arrow, anchor=west] (8.8, 0) -- node[above] {} (9.3, 0.5);
      \draw[arrow, anchor=west] (8.8, 0) -- node[below] {} (9.3, -0.5);
      \node[box, anchor=west] at (9.3, 0.7) {Holding costs\\ $\varphi^2: \mathbb{R}^{|\mathcal{H}|\times|\mathcal{P}|} \rightarrow \mathbb{R}$};
      \node[box, anchor=west] at (9.3, -0.7) {Stock-out costs\\ $\varphi^3: \mathbb{R}^{|\mathcal{H}|\times|\mathcal{P}|} \rightarrow \mathbb{R}$};
      \draw[arrow, anchor=west] (13.1, 0.5) -- node[above] {} (13.6, 0);
      \draw[arrow, anchor=west] (13.1, -0.5) -- node[below] {} (13.6, 0);
      \node[sum] (sum) at (13.6, 0) {$\Sigma$};
      \draw[arrow, anchor=west] (14, 0) -- node[above] {Prize} (16, 0);
      \draw[arrow, anchor=west] (14, 0) -- node[below] {$\theta_i$} (16, 0);
    \end{tikzpicture}
    \caption[Layers of our statistical model]{Layers of our statistical model.}
    \label{fig:statistical_model}
\end{figure}

To predict the prizes $\theta$, we design a \ac{PINN} $\varphi_w$, which maps each state $x^{\tau}$, including the set of historical observations $\hat{\mathcal{D}}^\tau_i$ for every customer $i \in \mathcal{V}_c$, to their corresponding prize values $\theta_i$. 
The \ac{PINN} leverages our understanding of the cost structure and effectively addresses data scarcity concerns.
For each quantile level $p\in \mathcal{P}$, we preprocess $\hat{\mathcal{D}}_i$ into demand quantiles given by
\begin{equation}
Q_p(\hat{\mathcal{D}}_i) = \inf \left\{d \in \hat{\mathcal{D}}_i : F(d) \geq p \right\}, 
\end{equation}
where $F(d)$ refers to the cumulative distribution function of the set $\hat{\mathcal{D}}_i$.
The first layer, denoted as~$\varphi^1$, acts as a demand estimation layer, tailoring the projection of future demand to the contextual information $\varPhi_i$ for look-ahead horizons $\mathcal{H}=\{0,..., H-1\}$. In our case, this corresponds to a linear combination of the individual features and their pairwise interactions as stated by
\begin{equation}
\varphi^1_{p}\left(\hat{\mathcal{D}}_i, \varPhi_i\right) = \text{ReLU}\left(Q_p(\hat{\mathcal{D}}_i) + w_1^{\top} \varPhi_i + \varPhi_i^{\top} w_2 \varPhi_i\right),
\end{equation}
where $w_1$ is a vector, and $w_2$ is a symmetric matrix with diagonal elements equal to zero.
Subsequently, two parallel layers come into play to estimate the holding \eqref{eq:holding} and stock-out costs \eqref{eq:stock_out}:
\begin{equation}
\label{eq:holding}
\varphi^2 = \sum_{p\in\mathcal{P}} {w_3^p}^{\top} \left(\text{ReLU}\left(I^{\tau}_i - \sum_{t=\tau}^{\tau+h} \varphi^1_{p}\left(\hat{\mathcal{D}}_i^\tau, \varPhi_i^{t}\right) \right)  \kappa_i\right)_{h\in\mathcal{H}}
\end{equation}
\begin{equation}\label{eq:stock_out}
\varphi^3 = \sum_{p\in\mathcal{P}} {w_4^p}^{\top} \left(\text{ReLU}\left(\sum_{t=\tau}^{\tau+h} \varphi^1_{p}\left(\hat{\mathcal{D}}_i^\tau, \varPhi_i^{t}\right) - I^{\tau}_i\right)  \kappa_i \rho\right)_{h\in\mathcal{H}}
\end{equation}
The resulting multi-layer statistical model $\varphi_{w} = \varphi^2 + \varphi^3$, with parameters $w=\{w_1\in\mathbb{R}^{|\varPhi_i|},w_2\in\mathbb{R}^{|\varPhi_i|\times|\varPhi_i|},w_3\in\mathbb{R}^{|\mathcal{H}|\times|\mathcal{P}|},w_4\in\mathbb{R}^{|\mathcal{H}|\times|\mathcal{P}|}\}$, is illustrated in \autoref{fig:statistical_model}.
This highly tailored design is justified by the inherent characteristics of the \ac{DSIRP}'s piecewise linear cost function, emerging from the distinct impacts of holding and stock-out costs, as inventory transitions from positive to negative.
The statistical model $\varphi_w$ can theoretically be replaced with any differentiable \ac{ML} model.
However, for the sake of interpretability, we make the choice to employ a statistical model where the core of layers $\varphi^1_{p},\varphi^2,\varphi^3$ is a generalized linear model offering streamlined training and involving a moderate number of interpretable parameters.

\section{Learning Algorithm}\label{sec:learning}
We now present our learning algorithm tailored to train policies for \ac{MDP}s with large state and action spaces $\calX$ and $\calU$ of the following form:
\begin{equation}\label{eq:CO_layer}
\delta_{w}(x^t) = \mathop{\arg \max}_{u \in \mathcal{U}(x^t)} {\theta}^\top g(u) + h(u), \quad \text{where} \quad \theta = \varphi_w(x^t).
\end{equation}
The learning algorithm is designed to select the parameters $w$ that result in an efficient policy $\delta_w : x_t \mapsto u_t$, where efficient is defined as minimizing the expected total (negative) reward.
\autoref{sub:setting} introduces an abstract \ac{MDP} setting, which clarifies the description of our learning algorithm, and highlights the main hypotheses needed to make it work.
\autoref{sub:imitated} introduces the imitated policy, and \autoref{sub:learningPb} formulates the learning problem. Next, \autoref{subsec:FYloss} defines a suitable loss function for training, and Sections~\ref{sub:Dagger}~and~\ref{sub:voting} discuss how to simultaneously train the model and include additional training samples with relevant states.

\subsection{Setting}
\label{sub:setting}

First, let's establish the generic setting to which our learning algorithm may be applied before delving into the specificities of our problem. For this purpose, we look at a  generic \ac{MDP} for which transitions \eqref{eq:deterministic_dynamic} and rewards \eqref{eq:deterministic_reward} are deterministic functions of the state~$x^t$, the decision $u^t$, and some exogenous noise $\xi^t$. 
\begin{equation}\label{eq:deterministic_dynamic}
    x^{t+1} = f(x^t,u^t,\xi^t)
\end{equation}
\begin{equation}\label{eq:deterministic_reward}
    r_t = \tilde r(x^t,u^t,\xi^t).
\end{equation}
By exogenous, we mean that $\xi^t$ is independent of the decisions taken. 
In other words, we consider an \ac{MDP} of the following form:
\begin{equation}\label{eq:genericMDP}
    \begin{aligned}
        \min_\delta \, &\bbE\Big[ \sum_{t=0}^{T-1}\tilde r\big(x^t,\delta(x^t),\xi^t \big) + \tilde r_T(x^T)\Big|x^0 \Big] \\
        \mathrm{s.t.}\, & \,\text{$u^t = \delta(x^t)$ and~\eqref{eq:deterministic_dynamic} $\forall t \in \{0,\dots,T\}$.}
    \end{aligned}
\end{equation}

An \emph{episode} then denotes a sequence $\bfxi = (\xi^0,\dots,\xi^{T-1})$ of noises. Given decisions $u^0,\dots,u^{T-1}$, we can compute the corresponding \emph{trajectory} $(x^0,u^0,\xi^0),\ldots, (x^{T-1},u^{T-1},\xi^{T-1})$ followed by the system.
To train our pipeline, we then require access to several episodes\:$\bfxi_1,\ldots,\bfxi_n$ of historical data, and an algorithm for the deterministic problem that arises from~\eqref{eq:genericMDP} when the distribution over $\bfxi$ is a Dirac, i.e., we have a single episode. 

In the case of our manuscript, this generic setting takes a specific form as outlined in \autoref{dsirp:sec:problemDescription}. Here, $\xi^t$ represents the uncertainty in demands $d(\xi^t)$, $u^t$ corresponds to routing decisions, and $x^t$ encapsulates the inventory, historical demand, and contextual information.The rewards~\eqref{eq:reward} and transitions~\eqref{eq:state_transition_contextual} are deterministic, and $\tilde r_T = 0$.

\subsection{Imitated Policy}
\label{sub:imitated}

Since our training set contains only episodes $\bfxi_i = (\xi_i^0,\dots,\xi_i^{T-1})$, we need to add decisions and turn these episodes into trajectories $(x_i^0,u_i^0,\xi_i^0),\ldots, (x_i^{T-1},u_i^{T-1},\xi_i^{T-1})$.
For that purpose, we determine the \emph{anticipative decisions} $\tilde \bfu^{\tau}(x^{\tau}; \bfxi) = \big(\tilde u^t(x^{\tau}; \bfxi)\big)_{t\in [\tau, T-1]}$. Given the current state $x^{\tau}$ at time $\tau$ and an episode $\bfxi$ for the rest of the time horizon, $\xi^{\tau},\dots,\xi^{T-1}$, it is the solution of the following deterministic problem.
\begin{equation}
    \label{eq:anticipative_decision}
    \begin{aligned}
        \tilde \bfu^{\tau}(x^{\tau}; \bfxi) \in \mathop{\argmin}_{u^\tau, ..., u^{T-1}}\, & \sum_{t=\tau}^{T-1} \tilde r(x^t,u^t,\xi^{t}) + \tilde r_T (x^T)\\
         \mathrm{s.t.}\,& \enskip x^{t} = f(x^{t-1},u^{t-1},\xi^{t-1}), \quad  \text{for all } t \in \{\tau+1,\dots,T\}
    \end{aligned} 
\end{equation}

A tailored algorithm is needed to solve this deterministic problem. In our case, this corresponds to a deterministic \ac{IRP}, which we solve with an adaptation of the exact algorithm of \citet{manousakis2021improved} that accommodates stock-outs. Any other available (heuristic or exact) algorithm could be used for this purpose. At time $\tau$, given the anticipative decisions $\tilde \bfu^{\tau}(x^{\tau}; \bfxi) = \big(\tilde u^t(x^{\tau}; \bfxi)\big)_{t\in [\tau, T-1]}$, our anticipative policy corresponds to the first time-step $\tilde u^\tau(x^{\tau}; \bfxi)$ as defined below.

\subsection{Learning Problem}
\label{sub:learningPb}

Let us now consider the anticipative policy, denoted by
$$ \delta^*_t (x^t; \bfxi) \mapsto  \tilde u^{t}(x^{t}; \bfxi).$$ 
This policy $\delta^*$ cannot be applied in practice since it requires information not available at time\:$t$. However, it can be recomputed a posteriori given an episode $\bfxi$.
Our learning algorithm trains our policy $\delta_w$, which uses only information available at time\:$t$, to imitate $\delta^*$.
In order to train a policy by imitation, we minimize a surrogate loss $\ell(u,\bar u)$ that penalizes taking a decision $u$ that does not correspond to the target decision $\bar u$. Given our pipeline, the decision taken by $\delta_w$ is a deterministic function. Therefore, we can consider a loss $\mathcal{L}(\theta,\bar u)$ as a function of $\theta$ and $\bar u$.
When we train a policy $\delta_w$ to imitate $\delta^*$, our objective is to find parameters $w$ that minimize an expected surrogate loss under their induced distribution of state
\begin{equation}\label{eq:imitation_learning_problem}
    \min_w \mathbb{E}_{\bfxi, X \sim d_w}\Big[\mathcal{L}\big(\varphi_w(X), \delta^*(X;\bfxi)\big)\Big],
\end{equation}
where $d_w$ is the distribution on $X$ induced by $w$, where we have aggregated the different time periods.
Practically, we do not know how to evaluate the expectation presented in~\eqref{eq:imitation_learning_problem}.
As such, we iteratively construct a training set $\mathcal{D}$ containing states $x$ sampled from $\delta_w$ and update the algorithm by solving the \emph{empirical imitation learning problem} given by
\begin{equation}
    \label{eq:empiricalLearningProblem}
    \min_w \sum_{(x,\bar u) \in \mathcal{D}} \mathcal{L}\big(\varphi_w(x),\bar u\big).
\end{equation}

\subsection{Fenchel-Young Loss}
\label{subsec:FYloss}

Next, we describe the surrogate loss $\mathcal{L}$ applied in our pipeline. 
Suppose we have a training set $(x^1,\bar u^1),\ldots, (x^n,\bar u^n)$ composed of states\:$x^k$ and their target decisions\:$\bar u^k$.
We want our policy~\eqref{eq:CO_layer} to output $\bar u^{k}$ when taking $x^k$ as input.
In other words, we want $\bar u^k$ to be an optimum of the \ac{CO}-layer.
Consequently, it is natural to use  as loss the non-optimality of $\bar u^k$ as a solution of the \ac{CO}-layer:
$$ \max_{u \in \mathcal{U}} \theta^\top g(u) + h(u) - \theta^\top  g(\bar u^k) - h(\bar u^k) \quad \text{where} \quad \theta = \varphi_w\left(x^k\right).$$
Based on this, we obtain a loss that is more robust to degeneracy when we perturb $\theta$ with a standard Gaussian random variable $Z$, leading to the Fenchel-Young loss\footnote{Omitting the Fenchel conjugate term.}
$$ \mathcal{L}(\theta,\bar u^k) = \mathbb{E} \Big[\max_{u \in \mathcal{U}} (\theta + Z)^\top g\left(u\right) + h\left(u\right) - \theta^\top  g(\bar u^k) - h(\bar u^k)\Big].$$
Several properties of this loss function are presented in \citet{baty2023combinatorial}. For our purposes, we need to underline that it is strongly convex and that its gradient is given by,
\begin{equation} 
\nabla_{\theta}\mathcal{L}(\theta,\bar u^k) = \mathbb{E}\big[g\left(o(\theta + Z)\right)\big] - g(\bar u^k),
\label{eq:FYLgradient}
\end{equation}
where $o$ is the \ac{CO}-layer's oracle that associates to $\theta$ an optimal solution of ${\max_{u \in \mathcal{U}} \theta^\top g(u) + h(u)}$.
Since we have such an oracle $o$ by hypothesis, we can compute a stochastic gradient by sampling on~$Z$.
Stochastic gradients in $w$ can, in this case, be computed by backpropagation~\eqref{eq:FYLgradient} using automatic differentiation.
The convexity of $\mathcal{L}$ ensures that the empirical imitation learning problem~\eqref{eq:empiricalLearningProblem} is convex when $\varphi_w$ is convex.

\subsection{Dataset Generation using DAgger Algorithm}
\label{sub:Dagger}

To generate the training dataset $\calD$, we use the DAgger algorithm proposed by \citet{ross2011reduction}. It is designed to ensure that the empirical learning problem~\eqref{eq:empiricalLearningProblem} provides a good solution to problem~\eqref{eq:imitation_learning_problem}.
The general idea of the algorithm is to draw trajectories according to a (mixture) policy given by
\begin{equation}\label{eq:mixture_policy}
\alpha \delta^* + (1-\alpha)\delta_w,
\end{equation}
where $\alpha$ is a mixture parameter.
Given an episode $\bfxi = (\xi^{0},\ldots,\xi^{T-1})$, a trajectory from policy~\eqref{eq:mixture_policy} can be sampled as follows. For $t$ going from\:$0$ to\:$T-1$, we select $\delta^*$ with probability $\alpha$. If so, we solve the deterministic problem~\eqref{eq:anticipative_decision} to take decision $u^t = \tilde u^{t}(x^{t}; \bfxi)$. Otherwise, we set $u^t = \delta_w(x^t)$. We then compute $x^{t+1}$ using the deterministic transition function~\eqref{eq:deterministic_dynamic}. 

\begin{algorithm}[H]
 	\caption{DAgger with anticipative policy} 
     \label{alg:DAgger_simple}
	\begin{algorithmic}[1]
    \State $\alpha_1,\ldots,\alpha_K$ a decreasing sequence of weights in $[0,1]$, $\mathcal{D} = \emptyset$, $w$ undefined.
		\For {$i=1,\dots, K$}
            \State Sample $\bfxi_i$, initial state $x^0$.
            \For {$t=0,\dots, T-1$}
                \State Transition to $x^{t+1}$ from policy~\eqref{eq:mixture_policy} with $\alpha_i$ and function\:\eqref{eq:deterministic_dynamic}.
                \State Add $(x^t,\tilde u^{t}(x^{t}; \bfxi_i))$ to $\mathcal{D}$.\label{op:agg_dataset}
            \EndFor
            \State Update the parameters $w$ solving learning problem~\eqref{eq:empiricalLearningProblem} with $\mathcal{D}$.
        \EndFor
	\end{algorithmic}
 \end{algorithm}
\autoref{alg:DAgger_simple}, referred to as \textit{Anticipative-DAgger}, starts with a collection of states generated by $\delta^*$ because $\delta_w$ is an inferior policy at the beginning of the training algorithm.
Iteration after iteration, $\delta_w$ improves, and we can therefore give more and more weight to $\delta_w$ in the mixture. In the last iterations, states come mostly from~$\delta_w$, which is aligned with the use of the induced distribution $d_w$ in \eqref{eq:imitation_learning_problem}.

An alternative would be the approach proposed in \citet{baty2023combinatorial}, where the data set of states and target decisions to imitate is defined a priori, i.e., at the beginning. Following \citet{baty2023combinatorial}'s learning paradigm, we sample initial states (with focus on inventories), calculating the anticipative decision \eqref{eq:anticipative_decision} for $\tau$ once and adding all pairs $(x^\tau,u^\tau),\dots,(x^{T-1},u^{T-1})$ to $\mathcal{D}$.
However, adding all pairs to the training set poses significant limitations
in our \ac{IRP} setting. Given the absence of an incentive to maintain inventory beyond period $T-1$, this triggers a reduction in delivered quantities during the final periods (as shown in \autoref{fig:eoh_effect} in our computational experiments), which biases the decisions.
Consequently, we modify the learning paradigm of \citet{baty2023combinatorial} and include only the first-period pair $(x^\tau,u^\tau)$ of state and anticipative decision to $\mathcal{D}$. 
In this approach, termed the \textit{Sampling} learning paradigm, states rely exclusively on the sampling strategy for initial states. This reliance may potentially lead to a lack of generalization.

Considering various look-ahead horizons between periods $1,\dots,T-1$ or even beyond $T-1$, to avoid the end-of-horizon effect, is theoretically possible.
However, this implies a trade-off between the end-of-horizon effect, the impact of the sampling strategy for the initial states, and the computational time for the anticipative decision. These considerations make determining the best and most efficient look-ahead horizon challenging.
The \textit{Anticipative-DAgger} learning paradigm works well in practice for our application.

\subsection{Voting Policy
}
\label{sub:voting}

Imitating the anticipative policy described in \autoref{sub:imitated} may seem surprising at first glance since such policies are known to be often sub-optimal.
An alternative would be to follow a voting policy, which consists in sampling several trajectories from any current state, solving them independently, and making them ``vote'' on the first decision by selecting the most frequent. Such an approach is generally more computationally demanding but leads to better policies. However, it cannot be applied to our setting due to the combinatorial size of the decision space; since the same decisions are highly unlikely to appear twice for different trajectories, a vote based on the highest frequency cannot be done. Nevertheless, to achieve a similar effect, we can replace step~\ref{op:agg_dataset} of~\autoref{alg:DAgger_simple} by the following steps, where $M$ is the number of trajectories used for the vote.
\begin{algorithm} \caption{DAgger with voting policy modifications}
	\begin{algorithmic}
    \State Sample $M$ trajectories $\tilde \bfxi_j$ for $j \in \{1,\ldots M\}$.
    \State Add $(x^t,\tilde u^{t}(x^{t}; \bfxi_j))$ to $\mathcal{D}$ for all $j \in \{1,\ldots M\}$.
	\end{algorithmic}
\end{algorithm}

The numerical experiments show that this improvement of Algorithm~\ref{alg:DAgger_simple} is computationally expensive during training but improves the performance of the policy learned.
Across all learning paradigms we employ a common \ac{ML} technique known as early stopping. Using validation episodes $\left(\bar{\bfxi}_j=(\bar{\xi}^0_j,\dots,\bar{\xi}^{T-1}_j)\right)_{j\in\{1,\dots,5\}}$ and an initial state $x^0$, we halt the training of our algorithm if there is no improvement in our current best parameters $w^*$ in terms of the total costs $\sum_{j=1}^{5}\sum_{t=0}^{T-1} \tilde r \big(x^t, \delta_{w^*}\left(x^t\right), \bar{\xi}^t_j\big)$, with transitions\:\eqref{eq:deterministic_dynamic}. In addition, for both \textit{DAgger-like} learning paradigms, we apply a reduction and subsampling strategy to the dataset. In general, we only consider the samples of the last ten epochs and furthermore select only 50\% of the state-decision pairs from the past epochs.

\section{Computational Experiments}
\label{sec:computational_exp}

All our experiments have been conducted using the \citet{gurobi} software as our \ac{MILP} solver for the \ac{CPCTSP}s\:\eqref{eq:CPCTSP} and the deterministic \ac{IRP}s needed to generate anticipative decisions\:\eqref{eq:anticipative_decision}. For training, we used a heterogeneous CPU cluster for increased computational capacity. Afterward, the trained pipeline and all the other benchmark methods were evaluated on a single thread of an AMD EPYC 7713P 64-core CPU. 

\subsection{Instance Design}
\label{sec:instances}

Instances from prior studies \citep[e.g.,][]{coelho2014heuristics} are limited in the characteristics of their demand scenarios. For a thorough experimental evaluation, we augmented these instances with additional demand scenarios and aligned their other parameters with those featured in other \ac{IRP} datasets \citep{archetti2007branch,archetti2011analysis,coelho2014heuristics}. These parameters include initial and maximum inventories, vehicle capacity, geographical locations, travel costs, and costs associated with inventory holding and stock-outs at customer locations. Furthermore, we assume an infinite production capacity and do not account for inventory holding costs at the supplier.

For an evaluation that permits even the most complex baseline methods (e.g., \ac{SAA}-based approaches) to produce solutions, we generated ten instances of each demand pattern with ten customer locations each. As explained in \autoref{subsec:Stat-Model}, we generated a set of 50 historical demand observations to enable initial estimates of demand quantiles. Moreover, we constructed two fixed sets of demand observations for each customer dedicated to validation and evaluation.
Our final evaluation covers a horizon spanning ten periods, enabling the \textsc{Saa-3} policy and anticipative baseline solutions to discover optimal solutions.

We also generate instances that have more varied demand characteristics than previous works focused on \textit{normal} and \textit{uniform} distribution. To that end, we incorporated instances featuring \textit{bimodal} demand distributions. Additionally, we generate instances with \textit{contextual} demand to gauge the effectiveness of our approach in handling contextual information.
For these contextual demand scenarios, we introduced features that influence demand, a level of control not typically available with real datasets where assumptions about the underlying contextual effect must be made.
All the specificities of these instances are thoroughly discussed in Appendix\:\ref{app:instances}.

\subsection{Benchmark Policies}

Many previous approaches are limited in their scope of applicability. Some methods assume a normal demand distribution \citep{huang2010modified}, whereas others are designed exclusively for discrete demand scenarios \citep{li2022two, bertazzi2013stochastic} or involve stock-out costs that are at least equivalent to the travel costs for a separate round trip from the depot to the customer \citep{gruler2018combining, bertazzi2013stochastic}.
Consensus strategies \citep{hvattum2009using, bent2004scenario}, as an alternative to the multi-scenario policies with their linkage constraints, become impractical when confronted with the extensive combinatorial nature of actions inherent in the \ac{SIRP}. Finally, in line with most \ac{IRP} definitions, the \ac{DSIRP} proposed by \citet{coelho2014heuristics} does not consider possible stock-outs when optimizing over the rolling horizon. However, from a cost perspective in a stochastic context, it is imperative to consider the possibility of failing to fulfill demand, especially for distant customers with low stock-out costs. 

Due to the above-mentioned limitations, our initial benchmark policies are rolling-horizon policies with a single demand episode inspired by \citet{coelho2014heuristics}.
The first method, the \textsc{Mean} policy, involves averaging all historical observations. The second method, denoted as \textsc{Saa-1}, directly selects a single historical observation. When contextual information is absent, we select the most recent observation. When contextual information is available, we select the observation with the smallest Euclidean distance based on the observed features. However, it is worth noting that these single-scenario methods have limitations: Averaging lacks the ability to consider contextual information, and relying on a single observation can lead to suboptimal decisions as it fails to capture the full range of demand uncertainty.
Therefore, we incorporated the approach proposed by \citet{bertazzi2013stochastic} into our set of baselines.
To tailor it to our needs, we have adapted their multi-scenario approach by using multiple demand samples for the current period, while approximating all subsequent periods with the mean but leveraging a \ac{MILP} formulation to solve it without the discrete demand assumption. This adaptation is called the \textsc{Saa-3} policy, drawing inspiration from the \ac{SAA} and employing three demand episodes.
As a pipeline-based benchmark policy, we use the approach proposed by \cite{baty2023combinatorial}.

\subsection{Operational Performance}
We first the operational performance of our pipeline. Subsequently, we discuss the training time here in \autoref{sec:training_time}.
Result\:\ref{res:performance_time} compares our pipeline, referred to as \textsc{ML-CO} policy and trained with the \textit{Voting-DAgger} learning paradigm, to \textsc{Mean, Saa-1, Saa-3} policy. Result\:\ref{res:EOH_effect}  confirms the end-of-horizon issue of \citet{baty2023combinatorial}'s learning paradigm.
Result\:\ref{res:DAgger_learning} compares the different learning paradigms. Result\:\ref{res:Voting_strategy} evaluates the voting policy within a \textit{DAgger-like} learning paradigm separately.
As the primary performance metric, we use the performance relative to the anticipatory baseline solution, evaluated as: $$\frac{\text{Total cost of the considered policy} - \text{Total cost of the anticipative policy}}{\text{Total cost of the anticipative policy}}.$$

\begin{figure}[!htb]
     \centering
     \input{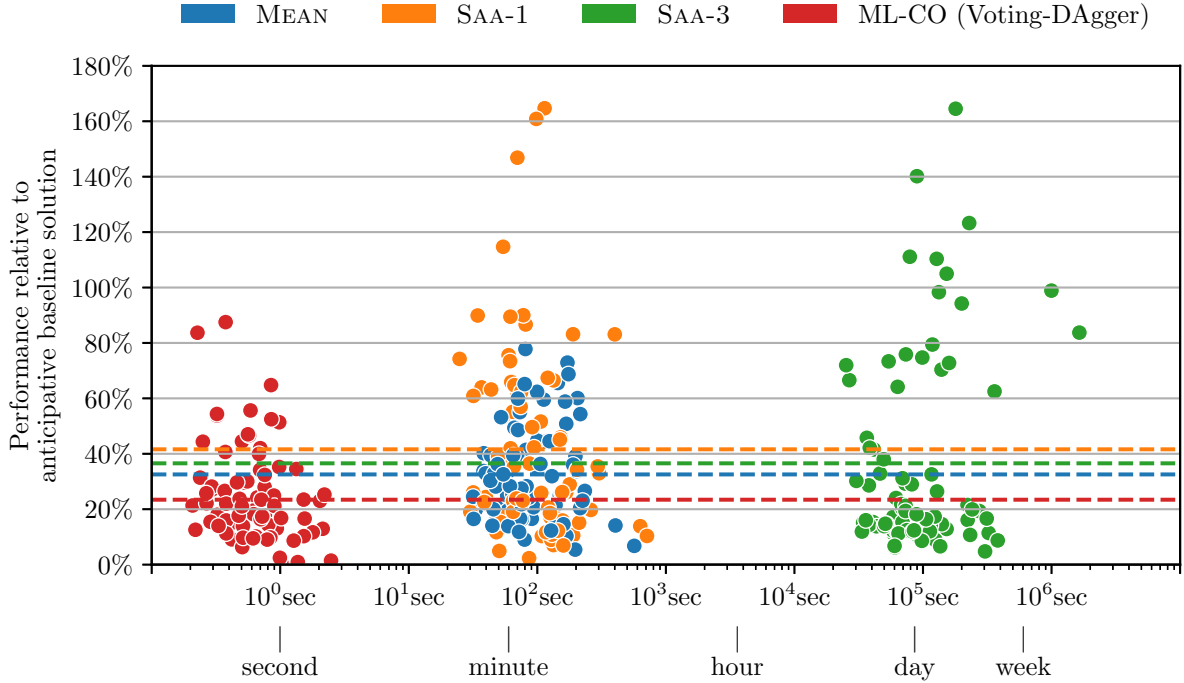}
     \caption[Performance and inference time per instance]{Performance and inference time per instance [horizontal line highlights the mean performance].}
     \label{fig:performance_overview}
 \end{figure}

\begin{result}\label{res:performance_time}
Our \textsc{ML-CO} policy outperforms all benchmarks on average across all penalties and demand patterns, with costs only 23.43\% higher compared to the anticipatory baseline solution. The best benchmark, the \textsc{Mean} policy, already incurs 32.55\% additional costs. In addition, inference times are an order of magnitude smaller than our benchmark policies. \autoref{fig:performance_overview} and \autoref{tab:result_1} summarize this major result.
\end{result}
 \begin{table}[!htb]
    \centering
        \caption[Average performance relative to anticipative baseline solution]{Average performance relative to anticipative baseline solution (in \%)}
    \label{tab:result_1}
    \begin{tabular}{llcccc}
\toprule
                                    &      &   \textsc{Mean} &  \textsc{Saa-1} &   \textsc{Saa-3} &  ML-CO \\
\textbf{demand pattern} & \textbf{shortage penalty} &        &        &         &             \\
\midrule
bimodal & low &  43.46 &  78.42 &   80.18 &       37.47 \\
                                    & high &  42.54 &  83.92 &  103.95 &       41.73 \\
\cline{3-6}
                                    &  &  43.00 &  81.17 &   92.07 &       39.60 \\
\midrule
contextual & low &  26.25 &  18.34 &   12.93 &        8.07 \\
                                    & high &  41.08 &  27.78 &   14.19 &        7.42 \\
\cline{3-6}
                                    &  &  33.66 &  23.06 &   13.56 &        7.74 \\
\midrule
normal & low &  17.51 &  16.22 &   11.91 &       13.82 \\
                                    & high &  17.95 &  18.39 &   14.78 &       18.65 \\
\cline{3-6}
                                    &  &  17.73 &  17.30 &   13.34 &       16.24 \\
\midrule
uniform & low &  27.58 &  34.54 &   27.98 &       24.51 \\
                                    & high &  44.07 &  55.30 &   26.56 &       35.76 \\
\cline{3-6}
                                    &  &  35.82 &  44.92 &   27.27 &       30.13 \\
\midrule
\multicolumn{2}{l}{\textit{average across scenarios}}  &     32.55 &  41.61 &   36.56 &       23.43 \\
\multicolumn{2}{l}{\textit{standard deviation across scenarios}}  &  (16.91) &  (34.48) &   (35.85) &       (18.84) \\
\bottomrule
\end{tabular}

\end{table}
As detailed in \autoref{tab:result_1}, none of the single scenario policies outperform the \textsc{ML-CO} policy on average across all penalties in any demand pattern.
For normal or uniform demand, the \textsc{Saa-3} policy 
is slightly better in terms of the additional costs incurred (16.24\% to 13.34\% for normal demand and 30.13\% to 27.27\% for uniform demand). However, bimodal demand or contextual information completely transforms the landscape (39.60\% to 92.07\% for bimodal demand and 7.74\% to 13.56\% for contextual demand).

\begin{figure}[!htb]
     \centering
     \input{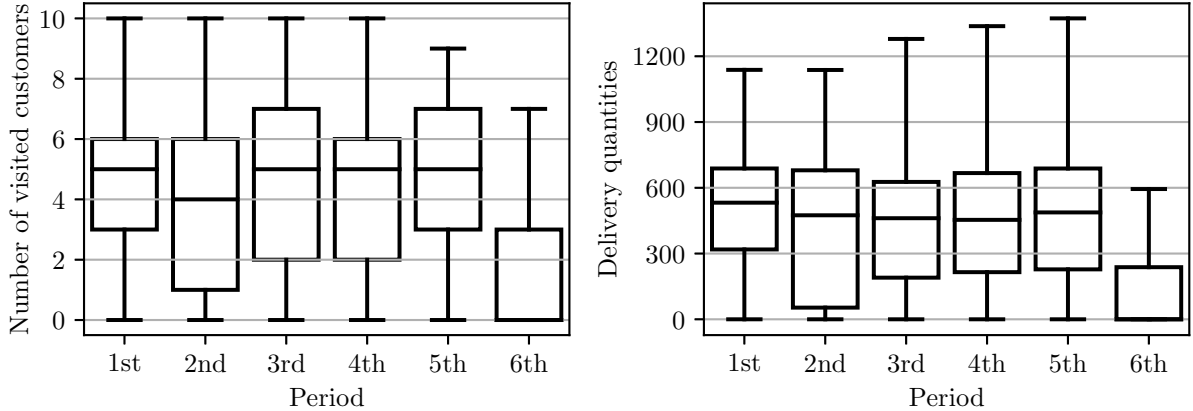}
     \caption{End-of-horizon effect using the learning paradigm of \citet{baty2023combinatorial}.}
     \label{fig:eoh_effect}
 \end{figure}

\begin{figure}[!htb]
     \centering
\input{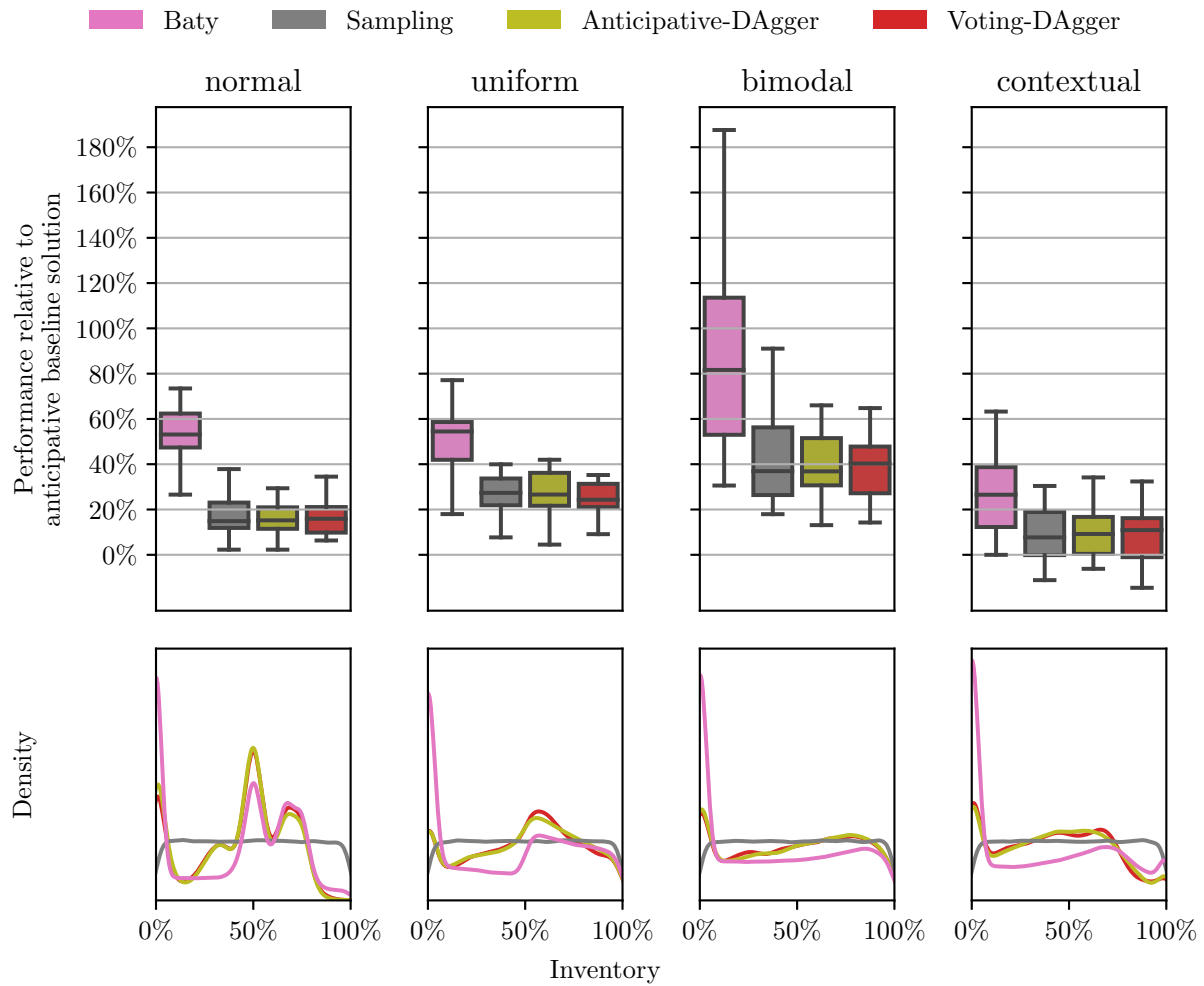}  
\caption[Performance per learning paradigm and demand pattern with its explored inventories]{Performance per learning paradigm and demand pattern with its explored inventories.}
     \label{fig:new_paradigm}
 \end{figure}

As depicted in \autoref{fig:performance_overview}, our pipeline has significant advantages regarding the inference times for dynamic decision-making. The inference times of our \textsc{ML-CO} policy are below one second. Even single-scenario policies like \textsc{Mean} and \textsc{Saa-1} require minutes to hours for inference, whereas the \textsc{Saa-3} policy requires inference times measured in days or weeks.
\begin{result}\label{res:EOH_effect}
The learning paradigm of \citet{baty2023combinatorial} suffers from the end-of-horizon effect. Given the absence of an incentive to maintain inventory beyond the final period in the anticipatory problem, the number of visited customers and the delivery quantities decrease in the later periods, as highlighted in \autoref{fig:eoh_effect}.
\end{result}
\begin{result}\label{res:DAgger_learning}
As summarized in \autoref{tab:performance_horizon_6} and \autoref{fig:new_paradigm}, both DAgger-based learning paradigms are superior. Regarding performance, \citet{baty2023combinatorial}'s learning paradigm is significantly worse. The \textit{Sampling} learning paradigm is comparable in performance but lacks robustness and generalizability of explored states.
\end{result}
\begin{table}[!htb]
    \centering
        \caption[Performance relative to anticipative baseline solution]{Performance relative to anticipative baseline solution (in \%)}
    \label{tab:performance_horizon_6}
\begin{tabular}{lccc}
\toprule
              &   & \textit{average}  &  \textit{standard} \\
\textbf{learning paradigm} & \textbf{look-ahead} &          &    \textit{deviation}                  \\
\midrule
Baty & 6 &    56.13 &               (34.38) \\
Sampling & 6 &    25.57 &               (23.21) \\
Anticipative-DAgger & 6 &    24.79 &               (19.73) \\
Voting-DAgger & 6 &    23.43 &               (18.84) \\
\toprule
\end{tabular}

\end{table}
As shown in \autoref{tab:performance_horizon_6}, the learning paradigm of \citet{baty2023combinatorial} performs worst (56.13\%), while the \textit{Sampling} learning paradigm (25.57\%) performs almost as good as the \textit{Anticipative-DAgger} (24.79\%)  and the \textit{Voting-DAgger} (23.43\%)\footnote{\label{note2}All learning paradigms use 600 samples in each epoch, except for the \textit{DAgger-like} learning paradigms, where the 600 samples are not reached until the 10th epoch due to the iterative dataset updates. If the voting policy is used with five scenarios per state, the number of states is reduced accordingly.}
\footnote{To demonstrate that the suboptimal performance of \citet{baty2023combinatorial}'s learning paradigm is unrelated to the difference between the look-ahead horizon (6 periods) and the evaluation horizon (10 periods), we include results with an evaluation horizon of 6 periods in Appendix\:\ref{app:results}. Expanding the look-ahead horizon is unfeasible due to computational constraints.}. But as depicted in \autoref{fig:new_paradigm}, the investigated states (inventories) of the \textit{Sampling} learning paradigm are based solely on the sampling strategy (here evenly distributed), lacking robustness and generalizability.

\begin{result}\label{res:Voting_strategy}
We improve the performance of the \textit{Anticipative-DAgger} learning paradigm from 24.79\% to 23.43\% by using the voting policy with five scenarios per state, as shown in \autoref{tab:performance_horizon_6}.
\end{result}
\begin{figure}[!htb]
     \centering
     \input{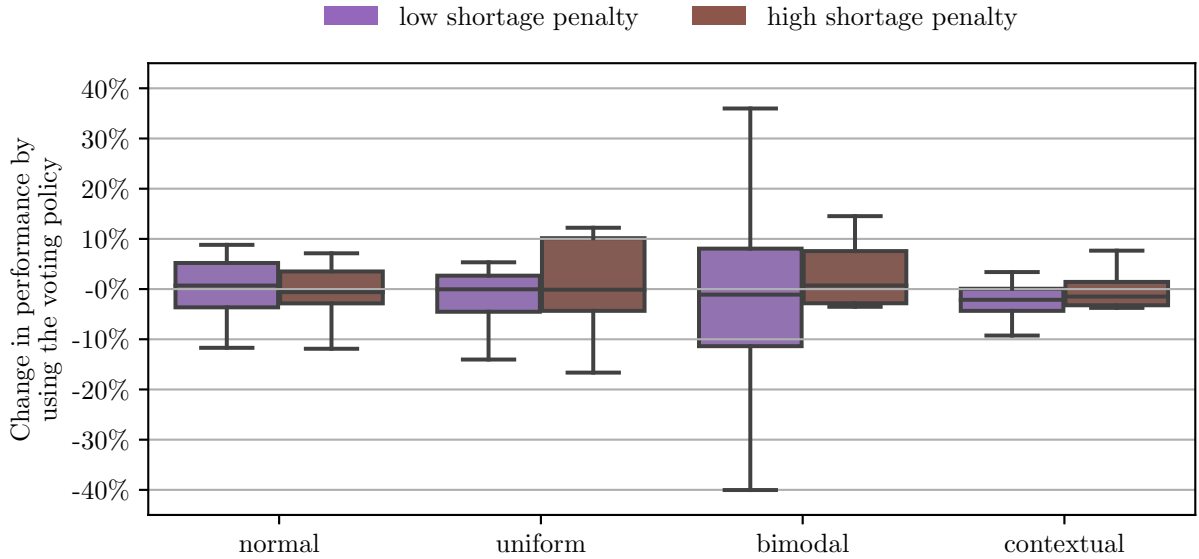}
     \caption{Impact of using \textit{Voting-DAgger} with five demand scenarios per state instead of \textit{Anticipative-DAgger}
     }
     \label{fig:voting_effect}
 \end{figure}
As detailed in \autoref{fig:voting_effect}, the voting policy with five scenarios instead of a single scenario per state is most valuable for demand with contextual information. No clear trend can be analyzed for the other demand patterns; this could be due to the reduced number of explored states using \textit{Voting-DAgger}.
 
\subsection{Training Time}\label{sec:training_time}
Training our pipeline demands substantial computational effort. Result\:\ref{res:training_time} examines the training time of the pipeline, which may span several days or weeks. To conclude, we analyze the effect of a reduced look-ahead horizon in Result\:\ref{res:reduced_look_ahead}.
\begin{result}\label{res:training_time}
As shown in \autoref{tab:time_horizon_6}, the learning paradigm of \citet{baty2023combinatorial} requires less time for sample generation (29.3 single-core CPU hours) compared to all other learning paradigms. This efficiency is achieved by extracting multiple samples from a single anticipatory decision problem. The time required for sample generation increases with the complexity of the learning paradigm and reaches 311.4 single-core CPU hours for \textit{Voting-DAgger}. In terms of policy update, both \textit{DAgger-like} learning paradigms require less time, as both use fewer samples up to the 10th epoch than the other learning paradigms.
 \begin{table}[!htb]
    \centering
        \caption{Average training time, normalized to a single CPU core (in hours).}
    \label{tab:time_horizon_6}
    \begin{tabular}{llrrr}
\toprule
              &   &  \textit{sample generation} &  \textit{policy update} &  \textit{others} \\
\textbf{learning paradigm} & \textbf{look-ahead} &                    &                &         \\
\midrule
Baty & 6 &               29.3 &         1157.7 &     2.3 \\
Sampling & 6 &              162.4 &         1285.8 &     2.8 \\
Anticipative-DAgger & 6 &              304.4 &          776.8 &     4.8 \\
Voting-DAgger & 6 &              311.4 &          996.2 &     6.9 \\
\bottomrule
\end{tabular}

\end{table}
\end{result}
Moreover, we observe that almost 100\% of the training time is spent on sample generation and policy updates, whereas the other steps take a time that is orders of magnitude smaller.

\begin{result}\label{res:reduced_look_ahead}
As demonstrated in \autoref{tab:time_horizon_3}, a reduced look-ahead horizon of three periods instead of six periods for the anticipatory decision problem drastically accelerates sample generation from 311.4 single-core CPU hours to 10.2 single-core CPU hours. According to \autoref{tab:performance_horizon_3}, this also leads to a better average performance. As detailed in \autoref{fig:result_4}, consolidating routing costs over time becomes less important with high shortage penalties, and demand without contextual information. However, performance deteriorates for contextual demand patterns.
  \begin{table}[!htb]
    \centering
        \caption{Average training time, normalized to a single CPU core (in hours).}
    \label{tab:time_horizon_3}
    \begin{tabular}{lcrrr}
\toprule
              &   &  \textit{sample generation} &  \textit{policy update} &  \textit{others} \\
\textbf{learning paradigm} & \textbf{look-ahead} &                    &                &         \\
\midrule
Voting-DAgger & 3 &               10.2 &          802.5 &     3.3 \\
\cmidrule{2-5}
              & 6 &              311.4 &          996.2 &     6.9 \\
\bottomrule
\end{tabular}

\end{table}
\begin{table}[!htb]
    \centering
        \caption[Performance relative to anticipative baseline solution]{Performance relative to anticipative baseline solution (in \%)}\label{tab:performance_horizon_3}
    \begin{tabular}{lccc}
\toprule
              &   &  \textit{average} &  \textit{standard} \\
\textbf{learning paradigm} & \textbf{look-ahead} &          &   \textit{deviation}                  \\
\midrule
Voting-DAgger & 3 &    23.29 &               (14.86) \\
\cline{2-4}
              & 6 &    23.43 &               (18.84) \\
\bottomrule
\end{tabular}

\end{table}
\begin{figure}[!htb]
     \centering
     \input{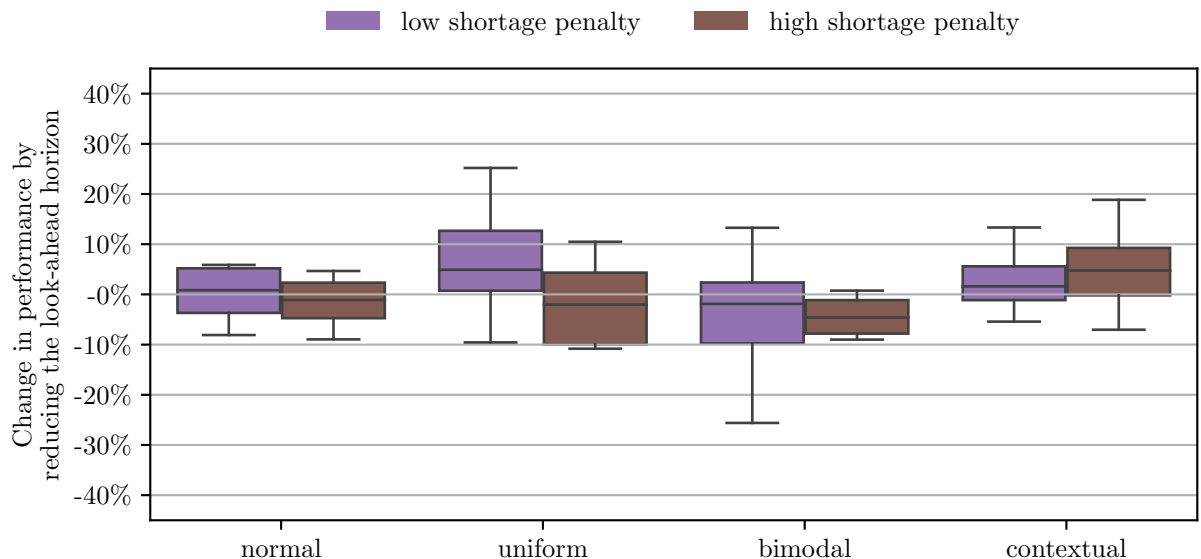}
     \caption[Impact of reducing the look-ahead horizon]{Impact of reducing the look-ahead horizon from 6 to 3 periods.}
     \label{fig:result_4}
 \end{figure}
\end{result}
\citet{coelho2014heuristics} report that the performance of their algorithm improves with a reduced look-ahead horizon. With our \textit{Voting-DAgger} learning paradigm, we could also observe this for non-contextual demand. However, for contextual demand, which was not studied by \citet{coelho2014heuristics}, we observe a performance degradation.
When setting up the pipeline, it is important to test different look-ahead horizons because the definition of the look-ahead horizon should be related to the relationship between demand and maximum inventory capacity.

\section{Conclusions}

Interest in dynamic decision-making has been rapidly growing, prompting ongoing efforts to develop more efficient solvers and heuristics. Our contribution to this landscape involves the introduction of a learning-based policy designed for the \ac{DSIRP}.
Compared to the benchmarks, our \textsc{ML-CO} policy achieved state-of-the-art performance and demonstrated a significant advantage in inference times for real-time decision-making. This comes, however, at the price of a longer training process, which must be performed in advance. We have identified two primary bottlenecks: the \ac{IRP} solver responsible for sample generation and the \ac{CPCTSP} oracle employed during policy updates. Replacing these computationally intensive \ac{MILP}s with efficient heuristics for the oracle, as done in \cite{baty2023combinatorial}, could significantly accelerate training.
Moreover, we suggest considering using \ac{ML} techniques such as stochastic gradient descent or transfer learning in future research. This implies that advancements in both research directions, \ac{CO} and \ac{ML}, have the potential to enhance the effectiveness of our approach.

Future work could consider more complex demand patterns or instance parameters. Although we have trained and evaluated our pipeline on instances with the same number of customers, it could be applied to networks of varying sizes.
A promising strategy could be to iteratively sample smaller subnetworks from large customer networks with thousands of customers, thereby creating a generalized model applicable to the entire network. The approach could also be extended to maximum-level policies, multi-depot settings, and multi-vehicle scenarios. Finally, beyond its application to the \ac{DSIRP}, our learning algorithm can train any policy encoded by a neural network with a \ac{CO} layer of the form~\eqref{eq:CO_layer}. Indeed, it suffices to be in the setting described in Section~\ref{sub:setting} for our algorithm to be applicable. Future works could explore the efficiency of our approach on other problems modeled as \ac{MDP}s with large state and action spaces. In conclusion, our approach achieves good performance for \ac{DSIRP}s and is versatile, paving the way for numerous promising avenues in future research.

\ACKNOWLEDGMENT{\noindent This work was supported by a fellowship within the IFI program of the German Academic Exchange Service (DAAD).
}

%
%
%

\renewcommand{\theHsection}{A\arabic{section}}
\begin{APPENDICES}
\section{Stochastic Instances}\label{app:instances}
We generated stochastic instances with \textit{normal}, \textit{uniform}, \textit{bimodal} and \textit{contextual} demand pattern, details are described in the following.
\paragraph{Normal and Uniform Demand}
Our stochastic instances with either \textit{normal} or \textit{uniform demand} for customer $i \in \mathcal{V}_c$ have been generated according to the following rules:
\begin{itemize}
    \item Mean demand $\mu_i$ is randomly selected from the set $\{10, \dots, 100\}$.
    \item Standard deviation $\sigma_i$ is randomly chosen from the set $\{2, \dots, 10\}$.
    \item Inventory capacity $C_i$ is determined by  $\mu_i$ times a random number from the set $\{2, 3, 4\}$.
    \item Initial inventory level $I_i^0$ is set as $C_i - \mu_i$.
    \item Per-unit holding cost $\kappa_i$ is randomly generated from a continuous uniform distribution within the interval $[0.02, 0.10]$.
    \item Shortage penalty $\rho$ determines the stock-out costs per unit, obtained by multiplying the holding cost $\kappa_i$ by the shortage penalty $\rho$. We evaluate two different scenarios with multiples of (i) $\rho=200$ (low shortage penalty) and (ii) $\rho=400$ (high shortage penalty).
    \item Vehicle capacity $B$ is calculated as 1.5 times the sum of the mean demands:
        \begin{align*}
        B = 1.5 \sum_{i\in\mathcal{V}_c} \mu_i
        \end{align*}
        This ensures that the vehicle capacity is appropriately scaled to accommodate the aggregated demand of all customers.
    \item Routing costs $\gamma_{ij}$, representing distances between vertices $i$ and $j$ have been calculated using the Euclidean distance. The coordinates of each vertex $i\in \mathcal{V}$ have been randomly sampled from a discrete uniform distribution within the interval $[0, 500]$.
    \item \textit{Normally-distributed demands} have been generated using a truncated normal distribution with mean $\mu_i$, standard deviation $\sigma_i$, and bounds of 0 and $C_i$.
    \item \textit{Uniformly-distributed demands} have been generated using a uniform distribution within the interval $\left[0, \frac{C_i}{2}\right]$. 
\end{itemize}

\paragraph{Bimodal Demand}
Our stochastic instances with a \textit{bimodal demand distribution} for each customer $i \in \mathcal{V}_c$ have been generated according to the following rules. The two components of the general mixture model have been generated using a truncated normal distribution with mean $\mu_i^1$, standard deviation $\sigma_i^1$ or mean $\mu_i^2$, standard deviation $\sigma_i^2$, and bounds of 0 and $C_i$.
\begin{itemize}
    \item The mean difference $\mu_{\Delta}$ has been randomly selected from the set $\{4, \dots, 20\}$.
    \item The mean $\mu_i^1$ has been randomly selected from the set $\{10, \dots, 50 - \mu_{\Delta}\}$.
    \item The mean $\mu_i^2$ has been randomly selected from the set $\{50 - \mu_{\Delta}, \dots, 100\}$.
    \item Standard deviations $\sigma_i^1$ and $\sigma_i^2$ were randomly chosen from the set $\{2, \dots, 10\}$.
    \item Inventory capacity $C_i$ have been set as $\frac{\mu_i^1+\mu_i^2}{2}$ times a random number from the set $\{2, 3, 4\}$. 
    \item The vehicle capacity $B$ has been set to 1.5 times the sum of the mean demands:
        \begin{align*}
        B = 1.5 \sum_{i\in\mathcal{V}_c} \frac{\mu_i^1+\mu_i^2}{2}
        \end{align*}
\end{itemize}
All other parameters were determined similarly to those for normal and uniform demand.

\paragraph{Contextual Demand}
We designed a synthetic dataset with eight distinct features, drawing inspiration from the generation process outlined by \citet{scikit-learn} in the $sklearn.datasets.make\_regression()$ method. This generation method is widely adopted in the scientific literature, as evidenced by references such as \cite{dessureault2022dpdr} and \cite{reddy1994using}.
To ensure diversity in our dataset, we first randomly selected a number of informative features for each instance from the set $\{2, \dots, 6\}$. This ensures that every instance includes at least two informative features, while the others remain non-informative.

We generated unscaled feature values within the range $[-1, +1]$. To diversify our instances, we employed one of three distributions (\textit{arcsin}, \textit{uniform}, or \textit{truncated normal}) with an equal likelihood of selection for each distribution. These distributions exert varying influence on the data's tails and centers. These features are then scaled by a randomly chosen factor from the set $\{\sqrt{10}, \dots, \sqrt{100}\}$, resulting in features $\lambda^1$ through $\lambda^8$ for each demand realization and each customer $i \ in \mathcal{V}_c$.
Furthermore, we introduced an exogenous noise $\xi$ generated in a manner consistent with the standard deviation of steady-state instances. Specifically, it follows a normal distribution with a mean of zero and a standard deviation randomly drawn from the set $\{2, \dots, 10\}$.
Let $M=\left\{(n, m) : n \in\{1, \dots, 8\}, m \in \{n+1,\dots, 8\}\right\}$ represent the set of all unordered pairs of features.
We constructed each \textit{contextual demand} realization for each customer $i$ as
\begin{equation}\label{eq:contextual_demand}
\min\left\{\max\left\{\mu_i + \sum_{n=1}^8 \alpha_n\lambda^n + \sum_{(n, m)\in M} \alpha_{nm} \lambda^{n} \lambda^{m} + \xi, 0\right\}, C_i\right\}.
\end{equation}
The coefficients $\alpha_{n}$ and $\alpha_{n,m}$ have been uniformly sampled once per instance in the interval $[-1, +1]$. These coefficients are the same among customers.
When dealing with the product of features, as in \autoref{eq:contextual_demand}, there is an equal likelihood of them being either informative or non-informative. In the latter case, we set $\alpha_{n,m} = 0$. As visible in the equation, each demand is bounded within the range of $0$ to $C_i$.
\section{Additional Results}\label{app:results}
\begin{figure}[!htb]
     \centering
     \input{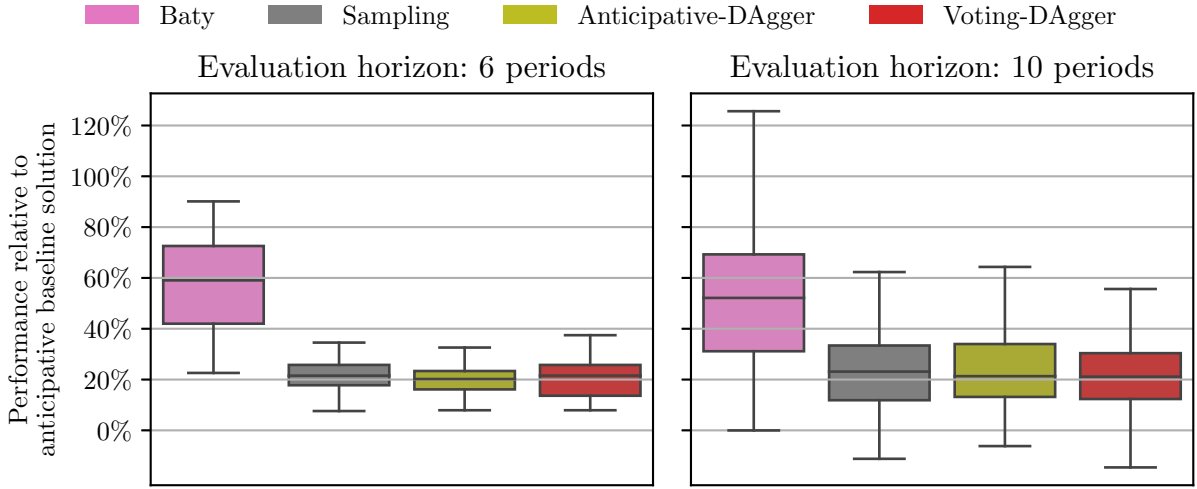}  \caption[Performance per learning paradigm and evaluation horizon]{Performance per learning paradigm and evaluation horizon [for all demand patterns].}
     \label{fig:new_paradigm_6}
 \end{figure}
\end{APPENDICES}


\bibliographystyle{informs2014trsc} 
\bibliography{literatur.bib} 


\end{document}